\documentclass[12pt]{article}
\usepackage{amsthm, amssymb, amsmath, latexsym}
\usepackage[utf8]{inputenc}
\usepackage[english]{babel}
\usepackage{longtable}

\begin{document}
\pagestyle{myheadings}
\thispagestyle{empty}
\setcounter{page}{1}

\newtheorem{definition}{Definition}
\newtheorem{proposition}{Proposition}
\newtheorem{theorem}{Theorem}
\newtheorem{lemma}{Lemma}
\newtheorem{corollary}{Corollary}
\newtheorem{remark}{Remark}
\newtheorem{algorithm}{Algorithm}
\theoremstyle{plain}
\mathsurround 2pt

\gdef\Aut{\mathop{\rm Aut}\nolimits}
\gdef\End{\mathop{\rm End}\nolimits}
\gdef\Ker{\mathop{\rm Ker}\nolimits}
\gdef\Im{\mathop{\rm Im}\nolimits}
\gdef\End{\mathop{\rm End}\nolimits}
\gdef\Inn{\mathop{\rm Inn}\nolimits}
\gdef\exp{\mathop{\rm exp}\nolimits}

\begin{center}
\textbf{\Large Groups of the nilpotency class $3$ of order $p^4$ as additive groups of local nearrings}
\end{center}

\begin{center}
\emph{Iryna Raievska, Maryna Raievska}\\
University of Warsaw, Warsaw, Poland;\\
Institute of Mathematics of National Academy of Sciences of Ukraine,\\
Kyiv, Ukraine\\
raeirina@imath.kiev.ua, raemarina@imath.kiev.ua
\end{center}

\date{}




\begin{abstract}
We consider groups of the nilpotency class $3$ of order $p^4$ which are the additive groups of local nearrings. It was shown that, for $p>3$, there exist a local nearring on one of such 4 groups.
\end{abstract}

\section{Introduction}

\

The classification of all nearrings up to certain orders is an open problem. It requires extensive computations, and the most suitable platform for their implementation is the computational algebra system GAP. Until well into 1990th the interest of pure mathematicians in nearring theory was stirred by, but in most cases also confined to the information that was produced by theoretical results for some research problems. However $27$ years ago the developers of the GAP~\cite{GAP} package SONATA~\cite{SONATA} shown the implementation of nearring theoretical algorithms. These are gradually becoming accepted both as standard tools for a working nearring theoretician, like certain methods of proof, and as worthwhile objects of study, like connections between notions expressed in theorems. The package SONATA provides methods for the construction and analysis of finite nearrings, as well as the library of all nearrings up to order $15$ and all nearrings with identity up to order $31$. The current version of the LocalNR package~\cite{LocalNR} (not yet redistributed with GAP) contains all local nearrings of order at most $361$, except those of orders $32$, $64$, $128$, $243$ and $256$. We have already calculated some classes of local nearrings of orders $32$, $64$, $128$, $243$ and $625$.

In \cite{IM_2021} it is proved that, up to isomorphism, there exist at least $p$ local nearrings on elementary abelian additive groups of order $p^3$, which are not nearfields. Lower bounds for the number of local nearrings on groups of order $p^3$ were obtained in~\cite{IMR_2022}. It is established that on each non-metacyclic non-abelian or metacyclic abelian groups of order $p^3$ there exist at least $p+1$ non-isomorphic local nearrings. It was proved that for $p>2$ every finite non-metacyclic $2$-generated $p$-group of nilpotency class $2$ with cyclic commutator subgroup is the additive group of a local nearring and in particular of a nearring with identity \cite{IM_20}.

Groups of the nilpotency class $2$ of order $p^4$ which are the additive group of local nearrings were investigated in~\cite{IMR_2023}. In this paper we consider groups of the nilpotency class $3$ of order $p^4$ which are the additive group of local nearrings. There are only examples of local nearrings on such groups constructed via GAP and the LocalNR package, i.e. of order 625 (see~\cite{Endom625}). It was shown that, for $p>3$, there exist a local nearring on one of such 4 groups.

\section{Preliminaries}

\

We will consider all groups of the nilpotency class $3$ of order $p^4$.

Let $[n,i]$ be the $i$-th group of order $n$ in the SmallGroups library in the computer system algebra GAP. We denote by $Z_n$ the cyclic group of order $n$.

It is an easy exercise for example in GAP to get the following assertions.

\begin{remark}\label{groups}
There are $3$ non-isomorphic groups of the nilpotency class $3$ of order $2^4=16$, which are:
\begin{enumerate}
  \item $D_{16}~[16,7]${\rm ;}
  \item $QD_{16}~[16,8]${\rm ;}
  \item $Q_{16}~[16,9]$.
\end{enumerate}
\end{remark}

\begin{remark}\label{groups}
There are $4$ non-isomorphic groups of the nilpotency class $3$ of order $3^4=81$, which are:
\begin{enumerate}
  \item $(C_3\times C_3\times C_3)\rtimes C_3~[81,7]${\rm ;}
  \item $(C_{9}\times C_3) \rtimes C_3~[81,8]${\rm ;}
  \item $(C_{9}\rtimes C_3)\rtimes C_3~[81,9]${\rm ;}
  \item $(C_{9}\rtimes C_3)\rtimes C_3~[81,10]$.
\end{enumerate}
\end{remark}

Using the packages SONATA and LocalNR it is possible to check the following assertion (see also~\cite{BN_98}).

\begin{theorem}
There is no local nearring on the groups of the nilpotency class $3$ of order of orders $2^4=16$ and $3^4=81$.
\end{theorem}

The following theorem contains the classification of groups of the nilpotency class $3$ of order $p^4$ with $p>3$ (see, for example, \cite{Burnside_1897}, \cite{Al-Hasanat_Almazaydeh_2022}).

\begin{theorem}\label{groups}
There are $4$ non-isomorphic groups of order $p^4$ with $p>3$, which are:
\begin{itemize}
  \item $H_1 = \langle a, b \colon a^p = b^p = c^p = [a, c]^p = [b, c] = e, [a, [a, c]] = [b, [a, c]] = e\rangle =(C_p\times C_p\times C_p)\rtimes C_p$, where $c = [a, b]${\rm ;}
  \item $H_2 = \langle a, b \colon a^{p^2}= b^p = [a, b]^p = [b, [a, b]] = e, [a, [a, b]] = a^p\rangle = (C_{p^2}\rtimes C_p)\rtimes C_p${\rm ;}.
  \item $H_3 = \langle a, b \colon a^{p^2}= b^p = [a, b]^p = [a, [a, b]] = e, [b, [a, b]] = a^p\rangle= (C_{p^2}\rtimes C_p)\rtimes C_p${\rm ;}
  \item $H_4 = \langle a, b \colon a^{p^2}= b^p = [a, b]^p = [a, [a, b]] = e, [b, [a, b]] = a^{2p}\rangle = (C_{p^2}\times C_p) \rtimes C_p$.
\end{itemize}
\end{theorem}

We will give the basic definitions.

\begin{definition}
A non-empty set $R$ with two binary operations $``+"$ and $``\cdot"$ is a \textbf{nearring} if:
\begin{description}
  \item[1)] $(R,+)$ is a group with neutral element $0${\rm ;}
  \item[2)] $(R,\cdot)$ is a semigroup{\rm ;}
  \item[3)] $x\cdot (y+z)=x\cdot y+x\cdot z$ for all $x$, $y$, $z\in R$.
\end{description}
Such a nearring is called a left nearring. If axiom 3) is replaced by an axiom $(x+y)\cdot z = x\cdot z + y\cdot z$ for all $x$, $y$, $z\in R$, then we get a right
nearring.
\label{nr}
\end{definition}

The group $(R,+)$ of a nearring $R$ is denoted by $R^+$ and called the {\em additive group} of $R$. It is easy to see that for each subgroup $M$ of $R^+$ and for each element $x\in R$ the set $xM=\{x\cdot y|y\in M\}$ is a subgroup of $R^+$ and in particular $x\cdot 0=0$. If in addition $0\cdot x=0$ for all $x\in R$, then the nearring $R$ is called {\em zero-symmetric}. Furthermore, $R$ is a {\em nearring with an identity} $i$  if the semigroup $(R,\cdot)$ is a monoid with identity element $i$. In the latter case the group of all invertible elements of the monoid $(R,\cdot)$ is denoted by $R^*$ and called the {\em multiplicative group} of $R$.  A subgroup $M$ of $R^+$ is called $R^*$-{\em invariant}, if $rM\leq M$ for each $r\in R^*$, and $(R,R)$-{\em subgroup}, if $xMy\subseteq M$ for arbitrary $x$, $y\in R$.\medskip

The following assertion is well-known (see, for instance, \cite{ClMal_66}, Theorem 3).

\begin{lemma}\label{exponent}
The exponent of the additive group of a finite nearring $R$ with identity $i$ is equal to the additive order of $i$ which coincides with the additive order of every invertible element of $R$.
\end{lemma}

\begin{definition}
A nearring $R$ with identity is called \textbf{local} if the set $L$ of all non-invertible elements of $R$ forms a subgroup of the additive group $R^{+}$.
\end{definition}

Throughout this paper $L$ will denote the subgroup of non-invertible elements of $R$.

The following lemma characterizes the main properties of finite local nearrings (see~\cite{AHS_2004}, Lemma~3.2).

\begin{lemma}\label{prop}
Let $R$ be a local nearring with identity $i$. Then the following statements hold{\rm :}
\begin{description}
\item[1)] $L$ is an $(R,R)$-subgroup of $R^{+}${\rm ;}
\item[2)] each proper $R^*$-invariant subgroup of $R^+$ is contained in $L${\rm ;}
\item[3)] the set $i+L$ forms a subgroup of the multiplicative group $R^*$.
\end{description}
\end{lemma}

Finite local nearrings with a cyclic subgroup of non-invertible elements are described in~\cite[Theorem~1]{RIM_2}.

\begin{theorem}\label{theorem_2}
Let $R$ be a local nearring of order $p^n$ with {$n\!> 1$} whose subgroup $L$ is cyclic and non-trivial. Then the additive group $R^+$ is either cyclic or is an elementary abelian group of order $p^2$. In the first case, $R$ is a commutative local ring, which is isomorphic to residual ring $\mathbb Z/p^n\mathbb Z$ with $n\ge 2$, in the other case there exist $p$ non-isomorphic such nearrings $R$ with $|L|=p$, from which $p-1$ are zero-symmetric nearrings and their multiplicative groups $R^{*}$ are isomorphic to a semidirect product of two cyclic subgroups of orders $p$ and $p-1$.
\end{theorem}

As a direct consequence of Theorem~\ref{theorem_2} we have the following result.

\begin{corollary}
Let $R$ be a local nearring of order $p^4$ with non-abelian additive group and is not a nearfield. Then the subgroup of non-invertible elements $L$ is a non-cyclic group of order $p^3$ or $p^2$.
\label{cor_1}
\end{corollary}

We define the binomial coefficient $\binom{n}{k}$ of integers $n$ and $k$ by

$\binom{n}{k}=\left\{
               \begin{array}{ll}
                 ${n!}/{k!(n-k)!}$, & if~\hbox{$0\leq k\leq n$;} \\
                    0, & if~\hbox{$k>n$ or $k<0$.}
                    \end{array}
                    \right.~\rm{(*)}$

\subsection{The group $H_1$}

\

Let $H_1$ be additively written group from Theorem~\ref{groups}. Then $H_1=\langle a\rangle +\langle b \rangle +\langle c \rangle +\langle d \rangle$for some elements $a$, $b$, $c$ and $d$ of $R$ satisfying the relations $ap=bp=cp=dp=0$, $b+c=c+b$, $a+d =d+a$, $b+d=d+b$, where $c=-a-b+a+b$ and $d=-a-c+a+c$.

\begin{lemma}\label{H1.1}
For arbitrary integers $k$ and $m$ in the group $H_1$ the equalities $-ak-cm+ak+cm=dkm$ and $cm+ak=ak+cm-dkm$ hold.
\end{lemma}

\begin{proof}
Since $-a-c+a+c=d$, we get $-c+a+c=a+d$. Then
$$-cm+ak+cm=(a+dm)k=ak+ckm.$$
Therefore, $-ak-cm+ak+cm=dkm$ and, so $cm+ak=ak+cm-dkm$.
\end{proof}

\begin{corollary}\label{H1.2}
For arbitrary integers $k$ and $m$ in the group $H_1$ the equalities $-cm+ak=ak-cm+dkm$ and $-cm-ak=-ak-cm-dkm$ hold.
\end{corollary}

\begin{lemma}\label{H1.3}
For arbitrary integers $k$ and $l$ in the group $H_1$ the equalities $-ak-bl+ak+bl=ckl+dl\binom{k}{2}$ and $bl+ak=ak+bl-ckl+dl\binom{k}{2}$ hold.
\end{lemma}

\begin{proof}
Let $k=1$. Since $-a-b+a+b=c$, we get $-b+a+b=a+c$. Then
$-bl+a+bl=a+cl=a+cl-dl\binom{1}{2}$, were $dl\binom{1}{2}=0$ due to ~\rm{(*)}.
Therefore, $-a-bl+a+bl=cl-dl\binom{1}{2}$ and, so $bl+a=a+bl-cl+dl\binom{1}{2}$.

The proof will be carried out by induction on $k$. For $k=1$ the equality is valid.
Let for $k$ the equality hold, i.e. $-bl+ak+bl=ak+clk-dl\binom{k}{2}$.

Let us prove the equality for $k+1$:
$$\begin{array}{l}
-bl+a(k+1)+bl=-bl+ak+a+bl=-bl+ak+bl+a+cl=\\
\qquad{}ak+clk-dl\binom{k}{2}+a+cl=ak+a+clk-dlk-dl\binom{k}{2}=\\
\qquad{}a(k+1)+cl(k+1)-dl(\binom{k}{2}+k)=a(k+1)+cl(k+1)-dl\binom{k+1}{2}.
\end{array}$$
\end{proof}

\begin{corollary}\label{H1.4}
For arbitrary integers $k$ and $l$ in the group $H_1$ the equalities $-bl+ak=ak-bl+ckl-dl\binom{k}{2}$ and $-bl-ak=-ak-bl-ckl-dl\binom{-k}{2}$ hold.
\end{corollary}

\begin{lemma}\label{H1.5}
For any natural numbers $k$, $l$, $n$, $m$ and $r$ in the group $H_1$ the following equality holds
$$\begin{array}{l}
(ak+bl+cm+dn)r=akr+blr+c(mr-kl\binom{r}{2})+d(nr+l\binom{k}{2}\binom{r}{2}-\\
\qquad{}km\binom{r}{2}+k^2l\binom{r}{r-3}).
\end{array}$$
\end{lemma}

\begin{proof}
The proof will be carried out by induction on $r$. For $r=1$ the equality is valid.
Let for $r$ the equality hold, i.e.
$$\begin{array}{l}
(ak+bl+cm+dn)r=akr+blr+c(mr-kl\binom{r}{2})+d(nr+l\binom{k}{2}\binom{r}{2}-\\
\qquad{}km\binom{r}{2}+k^2l\binom{r}{r-3}).
\end{array}$$
Using Lemmas~\ref{H1.1} and \ref{H1.3}, let us prove the equality for $r+1$:
$$\begin{array}{l}
(ak+bl+cm+dn)(r+1)=\\
\qquad{}akr+blr+c(mr-kl\binom{r}{2})+d(nr+l\binom{k}{2}\binom{r}{2}-\\
\qquad{}km\binom{r}{2}+k^2l\binom{r}{r-3})+ak+bl+cm+dn=\\
\qquad{}akr+blr+ak+c(m(r+1)-kl\binom{r}{2})-d(kmr-k^2l\binom{r}{2})+bl+\\
\qquad{}d(n(r+1)+l\binom{k}{2}\binom{r}{2}-km\binom{r}{2}+k^2l\binom{r}{r-3})=\\
\qquad{}akr+ak+blr+bl-cklr+dlr\binom{k}{2}+c(m(r+1)-kl\binom{r}{2})-\\
\qquad{}d(kmr-k^2l\binom{r}{2})+d(n(r+1)+l\binom{k}{2}\binom{r}{2}-km\binom{r}{2}+k^2l\binom{r}{r-3})=\\
\qquad{}a(k(r+1))+b(l(r+1))+c(m(r+1)-kl\binom{r+1}{2})+d(n(r+1)+\\
\qquad{}l\binom{k}{2}\binom{r+1}{2}-km\binom{r+1}{2}+k^2l\binom{r+1}{r-2}).
\end{array}$$
Therefore, the equality is valid for any $r$.
\end{proof}

\subsection{Nearrings with identity whose additive groups are isomorphic to $H_1$}

\

Let $R$ be a nearring with identity whose additive group of $R^+$ is isomorphic to $H_1$. Then $R^+=\langle a\rangle +\langle b \rangle +\langle c \rangle +\langle d \rangle$ for some elements $a$, $b$, $c$ and $d$ of $R$ satisfying the relations $ap=bp=cp=dp=0$, $b+c=c+b$, $a+d=d+a$,  $b+d=d+b$, where $c=-a-b+a+b$ and $d=-a-c+a+c$. In particular, each element $x\in R$ is uniquely written in the form $x=ax_1+bx_2+cx_3+dx_4$ with coefficients $0\le x_1<p$, $0\le x_2<p$, $0\le x_3<p$ and $0\le x_4<p$.

Since order of the element $a$ is equal to the exponent of group $G$, i.e. $p$, it follows that we can assume that $a$ is an identity of $R$, i.~e. $ax=xa=x$ for each $x\in R$. Furthermore, for each $x\in R$ there exist coefficients $\alpha(x)$, $\beta(x)$, $\gamma(x)$, $\phi(x)$  such that $xb=a\alpha(x)+b\beta(x)+c\gamma(x)+d\phi(x)$. It is clear that they are uniquely defined modulo $p$, so that some mappings $\alpha\colon R\to \mathbb Z_{p}$, ${\beta\colon R\rightarrow \mathbb Z_{p}}$, ${\gamma\colon R\rightarrow \mathbb Z_p}$, ${\phi\colon R\rightarrow \mathbb Z_p}$ are determined.

\begin{lemma}\label{H1.6}
Let $R$ be a nearring with identity whose additive group of $R^+$ is isomorphic to $H_1$. If $a$ coincides with identity element of $R$, $x=ax_1+bx_2+cx_3+dx_4$, $y=ay_1+by_2+cy_3+dy_4\in R$, $xb=a\alpha(x)+b\beta(x)+c\gamma(x)+d\phi(x)$, then
$$\begin{array}{l}
xy=a(x_1y_1+\alpha(x)y_2)+b(x_2y_1+\beta(x)y_2)+c(x_3y_1-x_1x_2\binom{y_1}{2}-x_2\alpha(x)y_1y_2+\\
\qquad{}\gamma(x)y_2)-\alpha(x)\beta(x)\binom{y_2}{2}+x_1\beta(x)y_3-x_2\alpha(x)y_3)+d(x_4y_1+\\
\qquad{}x_2\binom{x_1}{2}\binom{y_1}{2}-x_1x_3\binom{y_1}{2}+{x_1}^2x_2\binom{y_1}{y_1-3}+x_2y_1\binom{\alpha(x)y_2}{2}-\\
\qquad{}\alpha(x)x_3y_1y_2+x_1x_2\alpha(x)\binom{y_1}{2}y_2+\phi(x)y_2+\beta(x)\binom{\alpha(x)}{2}\binom{y_2}{2}-\alpha(x)\gamma(x)\binom{y_2}{2}+\\
\qquad{}\alpha(x)^2\beta(x)\binom{y_2}{y_2-3}+x_1\gamma(x)y_3-\beta(x)\binom{x_1}{2}y_3+\\
\qquad{}x_2\binom{\alpha(x)}{2}y_3-x_3\alpha(x)y_3+x_1^2\beta(x)y_4-x_1x_2\alpha(x)y_4).~\rm{(**)}
\end{array}$$
Moreover, for the mappings\\
$\begin{array}{l}
\alpha\colon R\to \mathbb \mathbb Z_{p}, {\beta\colon R\rightarrow \mathbb Z_{p}}, {\gamma\colon R\rightarrow \mathbb Z_p}, {\phi\colon R\rightarrow \mathbb Z_p}~\mbox{the following statements hold}{\rm :}\end{array}$
\begin{description}
  \item[\rm{(0)}] $\alpha(0)\equiv 0\; (\!\!\mod p)$, $\beta(0)\equiv 0\; (\!\!\mod p)$, ${\gamma(0)\equiv 0\; (\!\!\mod p)}${\rm ;}\\
   and ${\phi(0)\equiv 0\; (\!\!\mod p)}$ if and only if the nearring $R$ is zero-symmetric{\rm ;}
  \item[\rm{(1)}] $xc=c(x_1\beta(x)-x_2\alpha(x))+d(x_1\gamma(x)-\beta(x)\binom{x_1}{2}+x_2\binom{\alpha(x)}{2}-x_3\alpha(x))${\rm ;}
  \item[\rm{(2)}] $xd=d({x_1}^2\beta(x)-x_1x_2\alpha(x)){\rm ;}$
  \item[\rm{(3)}] $\alpha(xy)\equiv x_1\alpha(y)+\alpha(x)\beta(y)\; (\!\!\mod p\;){\rm ;}$
  \item[\rm{(4)}] $\beta(xy)\equiv x_2\alpha(y)+\beta(x)\beta(y)\; (\!\!\mod p)${\rm ;}
  \item[\rm{(5)}] $\gamma(xy)\equiv x_3\alpha(y)-x_1x_2\binom{\alpha(y)}{2}-x_2\alpha(x)\alpha(y)\beta(y)+\gamma(x)\beta(y)-\\
                   \alpha(x)\beta(x)\binom{\beta(y)}{2}+x_1\beta(x)\gamma(y)-x_2\alpha(x)\gamma(y)\; (\!\!\mod p\;){\rm ;}$
  \item[\rm{(6)}] $\phi(xy)\equiv x_4\alpha(y)+x_2\binom{x_1}{2}\binom{\alpha(y)}{2}-x_1x_3\binom{\alpha(y)}{2}+{x_1}^2x_2\binom{\alpha(y)}{\alpha(y)-3}+\\
x_2\alpha(y)\binom{\alpha(x)\beta(y)}{2}-\alpha(x)x_3\alpha(y)\beta(y)+x_1x_2\alpha(x)\binom{\alpha(y)}{2}y_2+\\
\phi(x)\beta(y)+\beta(x)\binom{\alpha(x)}{2}\binom{\beta(y)}{2}-\alpha(x)\gamma(x)\binom{\beta(y)}{2}+\\
\alpha(x)^2\beta(x)\binom{\beta(y)}{\beta(y)-3}+x_1\gamma(x)\gamma(y)-\beta(x)\binom{x_1}{2}\gamma(y)+\\
x_2\binom{\alpha(x)}{2}\gamma(y)-x_3\alpha(x)\gamma(y)+x_1^2\beta(x)\phi(y)-x_1x_2\alpha(x)\phi(y)\; (\!\!\mod p\;).$
\end{description}
\end{lemma}

\begin{proof}
Since $0\cdot a=a\cdot 0=0$, it follows that $R$ is a zero-symmetric nearring if and only if
$$0=0\cdot b=a\alpha(0)+b\beta(0)+c\gamma(0)+d\phi(0).$$
Equivalently we have
$$\alpha(0)\equiv 0\; (\!\!\mod p), {\beta(0)\equiv 0\; (\!\!\mod p)},$$
$${\gamma(0)\equiv 0\; (\!\!\mod p)}, {\phi(0)\equiv 0\; (\!\!\mod p)}.$$

Moreover, since $c=-a-b+a+b$, $d=-a-c+a+c$ and the left distributive law we have $0\cdot c=-0\cdot a-0\cdot b+0\cdot a+0\cdot b=0$ and $0\cdot d=-0\cdot a-0\cdot c+0\cdot a+0\cdot c=0$, whence $$0\cdot x=0\cdot (ax_1+bx_2+cx_3+dx_4)=(0\cdot a)x_1+(0\cdot b)x_2+(0\cdot c)x_3+(0\cdot d)x_4=0.$$  So that statement (0) holds.

Further, using Lemmas~\ref{H1.1} and \ref{H1.3}, we derive
$$\begin{array}{l}
xc=x(-a-b+a+b)=-xa-xb+xa+xb=-dx_4-cx_3-bx_2-ax_1-\\
\qquad{}d\phi(x)-c\gamma(x)-b\beta(x)-a\alpha(x)+ax_1+bx_2+cx_3+dx_4+\\
\qquad{}a\alpha(x)+b\beta(x)+c\gamma(x)+d\varphi(x)=-cx_3-bx_2-c\gamma(x)-\\
\qquad{}ax_1+dx_1\gamma(x)-b\beta(x)-a\alpha(x)+ax_1+bx_2+a\alpha(x)+\\
\qquad{}cx_3-dx_3\alpha(x)+b\beta(x)+c\gamma(x)=-cx_3-\\
\qquad{}bx_2-c\gamma(x)+dx_1\gamma(x)-b\beta(x)-ax_1+cx_1\beta(x)+d\beta(x)\binom{-x_1}{2}-\\
\qquad{}a\alpha(x)+a\alpha(x)+bx_2-cx_2\alpha(x)+dx_2\binom{\alpha(x)}{2}+\\
\qquad{}cx_3-dx_3\alpha(x)+b\beta(x)+c\gamma(x)=-cx_3-\\
\qquad{}bx_2-c\gamma(x)+dx_1\gamma(x)-b\beta(x)+cx_1\beta(x)-ax_1-\\
\qquad{}d{x_1}^2\beta(x)+d\beta(x)\binom{-x_1}{2}+ax_1+bx_2-cx_2\alpha(x)+\\
\qquad{}dx_2\binom{\alpha(x)}{2}+cx_3-dx_3\alpha(x)+b\beta(x)+c\gamma(x)=\\
\qquad{}c(x_1\beta(x)-x_2\alpha(x))+d(x_1\gamma(x)-\beta(x)\binom{x_1}{2}+x_2\binom{\alpha(x)}{2}-x_3\alpha(x)).
\end{array}$$

Using Lemmas~\ref{H1.1} and \ref{H1.3}, we get
$$\begin{array}{l}
xd=x(-a-c+a+c)=-xa-xc+xa+xc=-dx_4-cx_3-bx_2-ax_1-\\
\qquad{}d(x_1\gamma(x)-\beta(x)\binom{x_1}{2}+x_2\binom{\alpha(x)}{2}-x_3\alpha(x))-\\
\qquad{}c(x_1\beta(x)-x_2\alpha(x))+ax_1+bx_2+cx_3+dx_4+\\
\qquad{}c(x_1\beta(x)-x_2\alpha(x))+d(x_1\gamma(x)-\beta(x)\binom{x_1}{2}+x_2\binom{\alpha(x)}{2}-\\
\qquad{}x_3\alpha(x))=-cx_3-bx_2-c(x_1\beta(x)-x_2\alpha(x))-ax_1+d({x_1}^2\beta(x)-x_1x_2\alpha(x))+\\
\qquad{}ax_1+bx_2+cx_3+c(x_1\beta(x)-x_2\alpha(x))=d({x_1}^2\beta(x)-x_1x_2\alpha(x)).
\end{array}$$
We obtain statements~(1)--(2) of the lemma.

By the left distributive law and (1)--(2) of the lemma, we have
$$\begin{array}{l}
xy=x(ay_1)+x(by_2)+x(cy_3)+x(dy_4)=(ax_1+bx_2+cx_3+dx_4)y_1+\\
\qquad{}(a\alpha(x)+b\beta(x)+c\gamma(x)+d(\phi(x))y_2+(c(x_1\beta(x)-x_2\alpha(x))+d(x_1\gamma(x)-\\
\qquad{}\beta(x)\binom{x_1}{2})y_3+d(x_1^2\beta(x)-x_1x_2\alpha(x))y_4.
\end{array}$$

Furthermore, Lemma~\ref{H1.5} implies that
$$\begin{array}{l}(ax_1+bx_2+cx_3+dx_4)y_1=ax_1y_1+bx_2y_1+c(x_3y_1-x_1x_2\binom{y_1}{2})+\\
\qquad{}d(x_4y_1+x_2\binom{x_1}{2}\binom{y_1}{2}-x_1x_3\binom{y_1}{2}+\binom{y_1}{y_1-3}x_1^2x_2),
\end{array}$$

$$\begin{array}{l}(a\alpha(x)+b\beta(x)+c\gamma(x)+d\phi(x))y_2=a\alpha(x)y_2+b\beta(x)y_2+\\
\qquad{}c(\gamma(x)y_2+\alpha(x)\beta(x)\binom{y_2}{2})+d(\phi(x)y_2+\beta(x)\binom{\alpha(x)}{2}\binom{y_2}{2}-\\
\qquad{}\alpha(x)\gamma(x)\binom{y_2}{2}+\alpha(x)^2\beta(x)\binom{y_2}{y_2-3},
\end{array}$$

$$\begin{array}{l}(c(x_1b(x)-x_2\alpha(x))+d(x_1\gamma(x)-\beta(x)\binom{x_1}{2}+\\
\qquad{}x_2\binom{\alpha(x)}{2}-x_3\alpha(x)))y_3=c(x_1\beta(x)y_3-x_2\alpha(x)y_3)+\\
\qquad{}d(x_1\gamma(x)y_3-\beta(x)\binom{x_1}{2}y_3+x_2\binom{\alpha(x)}{2}y_3-x_3\alpha(x)y_3)
\end{array}$$

and

$$\begin{array}{l}d(x_1^2\beta(x)-x_1x_2\alpha(x))y_4=d(x_1^2\beta(x)y_4-x_1x_2\alpha(x)y_4).
\end{array}$$

By Lemmas~\ref{H1.1} and \ref{H1.3}, we have
$$\begin{array}{l}xy=a(x_1y_1+\alpha(x)y_2)+b(x_2y_1+\beta(x)y_2)+c(x_3y_1-x_1x_2\binom{y_1}{2}-x_2\alpha(x)y_1y_2+\\
\qquad{}\gamma(x)y_2-\alpha(x)\beta(x)\binom{y_2}{2}+x_1\beta(x)y_3-x_2\alpha(x)y_3)+d(x_4y_1+\\
\qquad{}x_2\binom{x_1}{2}\binom{y_1}{2}-x_1x_3\binom{y_1}{2}+{x_1}^2x_2\binom{y_1}{y_1-3}+x_2y_1\binom{\alpha(x)y_2}{2}-\\
\qquad{}\alpha(x)x_3y_1y_2+x_1x_2\alpha(x)\binom{y_1}{2}y_2+\phi(x)y_2+\beta(x)\binom{\alpha(x)}{2}\binom{y_2}{2}-\alpha(x)\gamma(x)\binom{y_2}{2}+\\
\qquad{}\alpha(x)^2\beta(x)\binom{y_2}{y_2-3}+x_1\gamma(x)y_3-\beta(x)\binom{x_1}{2}y_3+\\
\qquad{}x_2\binom{\alpha(x)}{2}y_3-x_3\alpha(x)y_3+x_1^2\beta(x)y_4-x_1x_2\alpha(x)y_4).
\end{array}$$

Finally, the associativity of multiplication in $R$ implies that
$$x(yb)=(xy)b=a\alpha(xy)+b\beta(xy)+c\gamma(xy)+d\phi(xy).$$
Furthermore, substituting $yb=a\alpha(y)+b\beta(y)+c\gamma(y)+d\phi(y)$ instead of $y$ in formula~$(**)$, we also have
$$\begin{array}{l}
x(yb)=a(x_1\alpha(y)+\alpha(x)\beta(y))+b(x_2\alpha(y)+\beta(x)\beta(y))+c(x_3\alpha(y)-\\
\qquad{}x_1x_2\binom{\alpha(y)}{2}-x_2\alpha(x)\alpha(y)\beta(y)+\gamma(x)\beta(y)-\\
\qquad{}\alpha(x)\beta(x)\binom{\beta(y)}{2}+x_1\beta(x)\gamma(y)-x_2\alpha(x)\gamma(y))+d(x_4\alpha(y)+\\
\qquad{}x_2\binom{x_1}{2}\binom{\alpha(y)}{2}-x_1x_3\binom{\alpha(y)}{2}+{x_1}^2x_2\binom{\alpha(y)}{\alpha(y)-3}+x_2\alpha(y)\binom{\alpha(x)\beta(y)}{2}-\\
\qquad{}\alpha(x)x_3\alpha(y)\beta(y)+x_1x_2\alpha(x)\binom{\alpha(y)}{2}y_2+\phi(x)\beta(y)+\beta(x)\binom{\alpha(x)}{2}\binom{\beta(y)}{2}-\\
\qquad{}\alpha(x)\gamma(x)\binom{\beta(y)}{2}+\alpha(x)^2\beta(x)\binom{\beta(y)}{\beta(y)-3}+x_1\gamma(x)\gamma(y)-\beta(x)\binom{x_1}{2}\gamma(y)+\\
\qquad{}x_2\binom{\alpha(x)}{2}\gamma(y)-x_3\alpha(x)\gamma(y)+x_1^2\beta(x)\phi(y)-x_1x_2\alpha(x)\phi(y)).
\end{array}$$
Comparing the coefficients under $a$, $b$, $c$ and $d$ in two expressions obtained for $x(yb)$, we derive statements~(3)--(6) of the lemma.
\end{proof}

\subsection{Local nearrings whose additive groups are isomorphic to $H_1$}

\

Let $R$ be a local nearring whose additive group of $R^+$ is isomorphic to $H_1$. Then $R^+=\langle a\rangle +\langle b \rangle +\langle c \rangle +\langle d \rangle$ for some elements $a$, $b$, $c$ and $d$ of $R$ satisfying the relations $ap=bp=cp=dp=0$, $b+c=c+b$, $a+d =d+a$  $b+d=d+b$, where $c=-a-b+a+b$ and $d=-a-c+a+c$. In particular, each element $x\in R$ is uniquely written in the form $x=ax_1+bx_2+cx_3+dx_4$ with coefficients $0\le x_1<p$, $0\le x_2<p$, $0\le x_3<p$ and $0\le x_4<p$.

Since order of the element $a$ is equal to the exponent of group $G$, i.e. $p$, it follows that we can assume that $a$ is an identity of $R$, i.~e. $ax=xa=x$ for each $x\in R$. Furthermore, for each $x\in R$ there exist coefficients $\alpha(x)$, $\beta(x)$, $\gamma(x)$, $\varphi(x)$  such that $xb=a\alpha(x)+b\beta(x)+c\gamma(x)+d\varphi(x)$. It is clear that they are uniquely defined modulo $p$, so that some mappings $\alpha\colon R\to \mathbb Z_{p}$, ${\beta\colon R\rightarrow \mathbb Z_{p}}$, ${\gamma\colon R\rightarrow \mathbb Z_p}$, ${\varphi\colon R\rightarrow \mathbb Z_p}$ are determined.

By Corollary~\ref{cor_1}, $L$ is the normal subgroup of order $p^3$ or $p^2$ in $R$.

Throughout this section let $R$ be a local nearring with $|R:L|=p$.

\begin{theorem}\label{H1.8}
Let $R$ be a local nearring whose additive group $R^{+}$ is isomorphic to a group $H_1$ and $|R:L|=p$. Then $R^{+}=\langle a\rangle +\langle b \rangle +\langle c \rangle +\langle d \rangle$, $ap=bp=cp=dp=0$, $b+c=c+b$, $a+d =d+a$  $b+d=d+b$, where $c=-a-b+a+b$, $d=-a-c+a+c$. If $a$ coincides with identity element of $R$, then the following statements hold:
\begin{itemize}
\item[1)] $L=\langle b\rangle + \langle c \rangle+ \langle d \rangle$ and $R^*=\{ax_1+bx_2+cx_3+cx_4\mid x_1\not\equiv 0 \; (\!\!\mod p\;)\}${\rm ;}
\item[2)] $xb=b\beta(x)+c\gamma(x)+d\varphi(x)${\rm ;}
\item[3)] $xc=cx_1\beta(x)+d(x_1\gamma(x)-\beta(x)\binom{x_1}{2})${\rm ;}
\item[4)] $xd=d{x_1}^2\beta(x)$.
\end{itemize}
\end{theorem}

\begin{proof}
Since $L$ consists the derived subgroup of $R^+$ it follows that the generators $b$, $d$ and $c$ we can choose such that $c=-a-b+a+b$ and $d=-a-c+a+c$. If $|L|=p^3$ then $L=\langle b\rangle + \langle c \rangle+ \langle d \rangle$ and $L$ is the $(R,R)$-subgroup in $R^+$ by statement~1) of Lemma~\ref{prop} it follows that $xb\in L$, hence $a\alpha(x)\in L$ for each $x\in R$. Thus $\alpha(x)\equiv 0 \; (\!\!\mod p)$, and so $xb=b\beta(x)+c\gamma(x)+d\phi(x)$ and statement 2) holds. Since $R^*=R\setminus L$  it follows ${R^*=\{ax_1+bx_2+cx_3+cx_4\mid x_1\not\equiv 0 \; (\!\!\mod p\;)\}}$ and $x=ax_1+bx_2+cx_3+dx_4$ is invertible if and only if $x_1\not\equiv 0 \; (\!\!\mod p\;)$, as claimed in statement 1).

Further, substituting $\alpha(x)\equiv 0 \; (\!\!\mod p)$ to 1) and 2) of Theorem~\ref{H1.6}, we derive $xc=cx_1\beta(x)+d(x_1\gamma(x)-\beta(x)\binom{x_1}{2})$ and $xd=d{x_1}^2\beta(x)$.
\end{proof}

As a consequence of Lemma~\ref{H1.6} and Theorem~\ref{H1.8} we have the following assertion.

\begin{corollary}\label{H1.9}
If $a$ coincides with identity element of $R$, $x=ax_1+bx_2+cx_3+dx_4$, $y=ay_1+by_2+cy_3+dy_4\in R$, $|R:L|=p$, then
$$\begin{array}{l}
xy=a(x_1y_1)+b(x_2y_1+\beta(x)y_2)+c(x_3y_1-x_1x_2\binom{y_1}{2}+x_1\beta(x)y_3+\\
\qquad{}\gamma(x)y_2)+d(x_4y_1+x_2\binom{x_1}{2}\binom{y_1}{2}-x_1x_3\binom{y_1}{2}+\\
\qquad{}{x_1}^2x_2\binom{y_1}{y_1-3}+\phi(x)y_2+x_1\gamma(x)y_3-\beta(x)\binom{x_1}{2}y_3+{x_1}^2\beta(x)y_4).~\rm{(***)}
\end{array}$$
Moreover, for the mappings $\beta\colon R\rightarrow \mathbb Z_p$, $\gamma\colon R\rightarrow \mathbb Z_p$, $\phi\colon R\rightarrow \mathbb Z_p$ the following statements hold{\rm :}
\begin{description}
  \item[\rm{(0)}] $\beta(0)\equiv 0\; (\!\!\mod p)$, ${\gamma(0)\equiv 0\; (\!\!\mod p)}$, ${\varphi(0)\equiv 0\; (\!\!\mod p)}$ if and only if the nearring $R$ is zero-symmetric;
  \item[\rm{(1)}] $\beta(xy)\equiv \beta(x)\beta(y)\; (\!\!\mod p)${\rm ;}
  \item[\rm{(2)}] $\gamma(xy)\equiv x_1\beta(x)\gamma(y)\; (\!\!\mod p)${\rm ;}
  \item[\rm{(3)}] $\phi(xy)\equiv \phi(x)\beta(y)+x_1\gamma(x)\gamma(y)-\beta(x)\binom{x_1}{2}\gamma(y)+x_1^2\beta(x)\phi(y)\; (\!\!\mod p)$.
\end{description}
\end{corollary}

Next, we give examples of local nearrings whose additive group $R^+$ is isomorphic to $H_1$.

\begin{lemma}\label{H1.10}
Let $R$ be a local nearring whose additive group of $R^+$ is isomorphic to $H_1$ and $|R:L|=p$. If $x=ax_1+bx_2+cx_3+dx_4$, $y=ay_1+by_2+cy_3+dy_4\in R$, then the mappings ${\beta\colon R\rightarrow \mathbb Z_p}$, ${\gamma\colon R\rightarrow \mathbb Z_p}$ and ${\phi\colon R\rightarrow \mathbb Z_p}$ from multiplication~$\rm{(***)}$ can be $\beta(x)\equiv {x_1}^2\; (\!\!\mod p\;)$ and
$\phi(x)=\left\{
           \begin{array}{ll}
             {x_2}^2, & if~\hbox{$x_1\equiv 0\; (\!\!\mod p\;)$}\\
             0, & if~\hbox{$x_1\not\equiv 0\; (\!\!\mod p\;)$}
           \end{array}
         \right.$, $\gamma(x)\equiv 0\; (\!\!\mod p\;)$.
\end{lemma}
\begin{proof}
It is easy to check that the functions $\beta$ and $\gamma$ satisfy conditions 1)--2) of Corollary~\ref{H1.9}. Since $\gamma(x)\equiv 0\; (\!\!\mod p\;)$ it follows that  $\phi(xy)=\phi(x)\beta(y)+x_1^2\beta(x)\phi(y)$. If $x_1y_1\not\equiv 0\; (\!\!\mod p\;)$, then $\phi(xy)=0=\phi(x)\beta(y)$. For $x_1y_1\equiv 0\; (\!\!\mod p\;)$ we need check three cases.

\begin{description}
  \item[1)] if $x_1\equiv 0\; (\!\!\mod p\;)$ and $y_1\not\equiv 0\; (\!\!\mod p\;)$, then $\phi(xy)=(x_2y_1)^2={x_2}^2{y_1}^2=\phi(x)\beta(y)${\rm ;}
  \item[2)] if $x_1\not\equiv 0\; (\!\!\mod p\;)$ and $y_1\equiv 0\; (\!\!\mod p\;)$, then $\phi(xy)=(x_2^2y_1)^2={x_2}^4{y_1}^3=\phi(x)\beta(y)${\rm ;}
  \item[3)] if $x_1\equiv 0\; (\!\!\mod p\;)$ and $y_1\equiv 0\; (\!\!\mod p\;)$, then $\phi(xy)=0=\phi(x)\beta(y)$.
\end{description}

We derive conditions 3) of Corollary~\ref{H1.9}.
\end{proof}

As a consequence of Lemma~\ref{H1.10} we have the following result.

\begin{theorem}\label{H1.11}
For each prime $p>3$ there exists a local nearring $R$ whose additive group $R^+$ is isomorphic to $H_1$.
\end{theorem}

\textbf{Example 1.} Let $G\cong (C_5\times C_5\times C_5)\rtimes C_5$. If $x=ax_1+bx_2+cx_3+dx_4$ and $y=ay_1+by_2+cy_3+dy_4\in G$ and $(G,+, \cdot)$ is a local nearring, then by Lemma~\ref{H1.10} $``\cdot"$ can be
$$\begin{array}{l}
x\cdot y=a(x_1y_1)+b(x_2y_1+\beta(x)y_2)+c(x_3y_1-x_1x_2\binom{y_1}{2}+x_1\beta(x)y_3)+\\
\qquad{}d(x_4y_1+x_2\binom{x_1}{2}\binom{y_1}{2}-x_1x_3\binom{y_1}{2}+\\
\qquad{}{x_1}^2x_2\binom{y_1}{y_1-3}+\phi(x)y_2-\beta(x)\binom{x_1}{2}y_3+{x_1}^2\beta(x)y_4),
\end{array}$$
where $\beta(x)\equiv {x_1}^2\; (\!\!\mod 5\;)$ and
$$\phi(x)=\left\{
           \begin{array}{ll}
             {x_2}^2, & if~\hbox{$x_1\equiv 0\; (\!\!\mod 5\;);$}\\
             0, & if~\hbox{$x_1\not\equiv 0\; (\!\!\mod 5\;)$.}
           \end{array}
         \right.$$

A computer program verified that for $p=5$ the nearring obtained in Lemma~\ref{H1.10} is indeed a local nearring (see Example~1), is deposited on GitHub:

\verb+https://github.com/raemarina/Examples/blob/main/LNR_625-7.txt+

From the package LocalNR and \cite{Endom625} we have the following number of all non-isomorphic zero-symmetric local nearrings on $H_1$ of order $625$.

\begin{center}
\begin{tabular}{|c|c|c|c|}
\hline
$IdGroup(R^+)$ & $IdGroup(R^*)$  & $StructureDescription(R^*)$                            & $n(R^*)$\\
\hline
$[625, 7]$     & $[500, 21]$     & $((C_5 \times C_5) \rtimes C_5) \rtimes C_4$           & 45  \\
103            & $[500, 42]$     & $C_5 \times ((C_5 \times C_5) \rtimes C_4)$            & 34 \\
               & $[500, 44]$     & $C_5 \times ((C_5 \times C_5) \rtimes C_4)$            & 6 \\
               & $[500, 46]$     & $(C_5 \times C_5 \times C_5) \rtimes C_4$              & 18 \\
\hline
\end{tabular}
\end{center}

\footnotesize{{\bf Acknowledgement.} The authors are grateful IIE-SRF for the support of their fellowship at the University of Warsaw.

\end{document}